\input amstex
\documentstyle{amsppt}
%
\catcode`@=11
\redefine\output@{%
  \def\break{\penalty-\@M}\let\par\endgraf
  \ifodd\pageno\global\hoffset=105pt\else\global\hoffset=8pt\fi  
  \shipout\vbox{%
    \ifplain@
      \let\makeheadline\relax \let\makefootline\relax
    \else
      \iffirstpage@ \global\firstpage@false
        \let\rightheadline\frheadline
        \let\leftheadline\flheadline
      \else
        \ifrunheads@ 
        \else \let\makeheadline\relax
        \fi
      \fi
    \fi
    \makeheadline \pagebody \makefootline}%
  \advancepageno \ifnum\outputpenalty>-\@MM\else\dosupereject\fi
}
\catcode`\@=\active
\nopagenumbers
\def\negskp{\hskip -2pt}
\def\MatGrU{\operatorname{U}}
\def\MatGrSU{\operatorname{SU}}
\def\vtrule{\vrule height 12pt depth 6pt}
\def\vtttrule{\vrule height 12pt depth 19pt}
\def\boxit#1#2{\vcenter{\hsize=122pt\offinterlineskip\hrule
  \line{\vtttrule\hss\vtop{\hsize=120pt\centerline{#1}\vskip 5pt
  \centerline{#2}}\hss\vtttrule}\hrule}}
\def\chirk{\special{em:point 1}\kern 1.2pt\raise 0.6pt
  \hbox to 0pt{\special{em:point 2}\hss}\kern -1.2pt
  \special{em:line 1,2,0.3pt}\ignorespaces}
\def\Chirk{\special{em:point 1}\kern 1.5pt\raise 0.6pt
  \hbox to 0pt{\special{em:point 2}\hss}\kern -1.5pt
  \special{em:line 1,2,0.3pt}\ignorespaces}
\def\Hirk{\kern 0pt\special{em:point 1}\kern 4pt\special{em:point 2}
  \kern -3pt\special{em:line 1,2,0.3pt}\ignorespaces}
\accentedsymbol\uuud{d\unskip\kern -3.8pt\raise 1pt\hbox to
  0pt{\chirk\hss}\kern 0.2pt\raise 1.9pt\hbox to 0pt{\chirk\hss}
  \kern 0.2pt\raise 2.8pt\hbox to 0pt{\chirk\hss}\kern 3.5pt}
\accentedsymbol\bolduuud{\bold d\unskip\kern -4pt\raise 1.1pt\hbox to
  0pt{\chirk\hss}\kern 0pt\raise 1.9pt\hbox to 0pt{\chirk\hss}
  \kern 0pt\raise 2.5pt\hbox to 0pt{\chirk\hss}\kern 4pt}
\accentedsymbol\uud{d\unskip\kern -3.72pt\raise 1.5pt\hbox to
  0pt{\chirk\hss}\kern 0.26pt\raise 2.4pt\hbox to 0pt{\chirk\hss}
  \kern 3.5pt}
\accentedsymbol\bolduud{\bold d\unskip\kern -4pt\raise 1.4pt\hbox to
  0pt{\chirk\hss}\kern 0pt\raise 2.3pt\hbox to 0pt{\chirk\hss}
  \kern 4pt}
\accentedsymbol\uuuD{D\unskip\kern -4.6pt\raise 2pt\hbox to
  0pt{\Chirk\hss}\kern 0.2pt\raise 3pt\hbox to 0pt{\Chirk\hss}
  \kern 0.2pt\raise 4.1pt\hbox to 0pt{\Chirk\hss}\kern 4.5pt}
\accentedsymbol\uuD{D\unskip\kern -4.5pt\raise 2.4pt\hbox to
  0pt{\Chirk\hss}\kern 0.2pt\raise 3.4pt\hbox to 0pt{\Chirk\hss}
  \kern 4.4pt}
\accentedsymbol\uD{D\unskip\kern -4.4pt\raise 3pt\hbox to
  0pt{\Chirk\hss}\kern 4.4pt}
\accentedsymbol\bolduuuD{\bold D\unskip\kern -4.6pt\raise 2pt\hbox to
  0pt{\Chirk\hss}\kern 0pt\raise 3pt\hbox to 0pt{\Chirk\hss}
  \kern 0pt\raise 4.1pt\hbox to 0pt{\Chirk\hss}\kern 4.5pt}
\accentedsymbol\bolduuD{\bold D\unskip\kern -4.5pt\raise 2.4pt\hbox to
  0pt{\Chirk\hss}\kern 0pt\raise 3.4pt\hbox to 0pt{\Chirk\hss}
  \kern 4.4pt}
\accentedsymbol\bolduD{\bold D\unskip\kern -4.5pt\raise 3pt\hbox to
  0pt{\Chirk\hss}\kern 4.4pt}
\accentedsymbol\uuuU{U\unskip\kern -5.1pt\raise 2pt\hbox to
  0pt{\Chirk\hss}\kern 0.2pt\raise 3pt\hbox to 0pt{\Chirk\hss}
  \kern 0.2pt\raise 4.1pt\hbox to 0pt{\Chirk\hss}\kern 4.9pt}
\accentedsymbol\uuU{U\unskip\kern -5.0pt\raise 2.4pt\hbox to
  0pt{\Chirk\hss}\kern 0.2pt\raise 3.4pt\hbox to 0pt{\Chirk\hss}
  \kern 4.9pt}
\accentedsymbol\uU{U\unskip\kern -4.9pt\raise 3pt\hbox to
  0pt{\Chirk\hss}\kern 4.8pt}
\accentedsymbol\bolduPsi{\boldsymbol\Psi\unskip\kern -6.5pt\raise 3pt
  \hbox to 0pt{\Hirk\hss}\kern 6.6pt}
\accentedsymbol\bolduuPsi{\boldsymbol\Psi\unskip\kern -6.5pt\raise 3pt
  \hbox to 0pt{\Hirk\hss}\unskip\kern 0pt\raise 3.8pt
  \hbox to 0pt{\Hirk\hss}\kern 6.6pt}
\accentedsymbol\bolduuuPsi{\boldsymbol\Psi\unskip\kern -6.5pt\raise 3pt
  \hbox to 0pt{\Hirk\hss}\unskip\kern 0pt\raise 3.8pt
  \hbox to 0pt{\Hirk\hss}\unskip\kern 0pt\raise 4.6pt
  \hbox to 0pt{\Hirk\hss}\kern 6.6pt}
\def\blue#1{#1}
\catcode`#=11\def\diez{#}\catcode`#=6
\catcode`_=11\def\podcherkivanie{_}\catcode`_=8
\def\mycite#1{\cite{\blue{#1}}\immediate\special{ps:
     ShrHPSdict begin /ShrBORDERthickness 0 def}}

\def\mytag#1{%
    \tag#1}
\def\mythetag#1{\thetag{\blue{#1}}\immediate\special{ps:
     ShrHPSdict begin /ShrBORDERthickness 0 def}}
\def\myrefno#1{\no#1}
\def\myhref#1#2{\blue{#2}\immediate\special{ps:
     ShrHPSdict begin /ShrBORDERthickness 0 def}}
\def\myEarXivlink{\myhref{http://arXiv.org}{http:/\negskp/arXiv.org}}
\def\myGeoCities{\myhref{http://www.geocities.com}{GeoCities}}
\def\mytheorem#1{\csname proclaim\endcsname{Theorem #1}}
\def\mythetheorem#1{\blue{#1}\immediate\special{ps:
     ShrHPSdict begin /ShrBORDERthickness 0 def}}
\def\mylemma#1{\csname proclaim\endcsname{Lemma #1}}

\def\mycorollary#1{\csname proclaim\endcsname{Corollary #1}}

\def\mydefinition#1{\definition{Definition #1}}
\def\mythedefinition#1{\blue{#1}\immediate\special{ps:
     ShrHPSdict begin /ShrBORDERthickness 0 def}}

\pagewidth{360pt}
\pageheight{606pt}
\topmatter
\title
The electro-weak and color bundles 
for the Standard Model in a gravitation field.
\endtitle
\author
R.~A.~Sharipov
\endauthor
\address 5 Rabochaya street, 450003 Ufa, Russia\newline
\vphantom{a}\kern 12pt Cell Phone: +7-(917)-476-93-48
\endaddress
\email \vtop to 30pt{\hsize=280pt\noindent
\myhref{mailto:r-sharipov\@mail.ru}
{r-sharipov\@mail.ru}\newline
\myhref{mailto:ra\podcherkivanie sharipov\@lycos.com}{ra\_\hskip 1pt
sharipov\@lycos.com}\newline
\myhref{mailto:R\podcherkivanie Sharipov\@ic.bashedu.ru}
{R\_\hskip 1pt Sharipov\@ic.bashedu.ru}\vss}
\endemail
\urladdr
\vtop to 20pt{\hsize=280pt\noindent
\myhref{http://www.geocities.com/r-sharipov}
{http:/\negskp/www.geocities.com/r-sharipov}\newline
\myhref{http://www.freetextbooks.boom.ru/index.html}
{http:/\negskp/www.freetextbooks.boom.ru/index.html}\vss}
\endurladdr
\abstract
    It is known that the Standard Model describing all of the currently
known elementary particles is based on the $\MatGrU(1)\times\MatGrSU(2)
\times\MatGrSU(3)$ symmetry. In order to implement this symmetry on the
ground of a non-flat space-time manifold one should introduce three
special bundles. Some aspects of the  mathematical theory of these bundles
are studied in this paper.
\endabstract
\endtopmatter
\loadbold
\loadeufb
\TagsOnRight
\document

\rightheadtext{The electro-weak and color bundles \dots}
\head
1. The Dirac bundle and its basic fields. 
\endhead
    Spin is a common property of all elementary particles. It 
is described in terms of the Dirac bundle $DM$, where $M$ is a
space-time manifold. The Dirac bundle is a four-dimensional complex 
bundle over the four-dimensional real manifold $M$. The detailed 
description of this bundle can be found in \mycite{1}, \mycite{2}, 
and \mycite{3}. In this section we shall remind in brief some of its
properties.\par
    The base manifold $M$ of the Dirac bundle $DM$ is equipped with 
the following three geometric structures which are well-known in 
General Relativity:
\roster
\item a pseudo-Euclidean Minkowski-type metric $\bold g$;
\item an orientation;
\item a polarization
\endroster
(see details in \mycite{4}). The Dirac bundle $DM$ is linked to the
above three structures through frame pairs. The most popular frame
pairs are listed in the following table:
$$
\hskip -2em
\aligned
&\boxit{Canonically orthonormal}{chiral frames}\to
\boxit{Positively polarized}{right orthonormal frames}\\
&\boxit{$P$-reverse}{anti-chiral frames}\to
\boxit{Positively polarized}{left orthonormal frames}\\
&\boxit{$T$-reverse}{anti-chiral frames}\to
\boxit{Negatively polarized}{right orthonormal frames}\\
&\boxit{$PT$-reverse}{chiral frames}\to
\boxit{Negatively polarized}{left orthonormal frames}
\endaligned
\mytag{1.1}
$$
(see details in \mycite{1} and \mycite{2}). Each of the frames
in the right hand side of the table \mythetag{1.1} is an ordered
set of four vector fields $\boldsymbol\Upsilon_0$, $\boldsymbol
\Upsilon_1$, $\boldsymbol\Upsilon_2$, $\boldsymbol\Upsilon_3$
linearly independent at each point $p$ of some open domain 
$U\subset M$. We denote it $(U,\,\boldsymbol\Upsilon_0,\,\boldsymbol
\Upsilon_1,\,\boldsymbol\Upsilon_2,\,\boldsymbol\Upsilon_3)$. Similarly,
each of the frames listed in the left hand side of the table
\mythetag{1.1} is an ordered set of four spinor fields $\boldsymbol\Psi_1$,
 $\boldsymbol\Psi_2$, $\boldsymbol\Psi_3$, $\boldsymbol\Psi_4$, i\.\,e\.
four smooth sections of the Dirac bundle $DM$, linearly independent at 
each point $p\in U$. We denote it $(U,\,\boldsymbol\Psi_1, \,\boldsymbol
\Psi_2,\,\boldsymbol\Psi_3,\,\boldsymbol\Psi_4)$. The spinor frames 
in the left hand side of the table \mythetag{1.1} are defined with the use
of the following basic spin-tensorial fields (see \mycite{2}):
$$
\vcenter{\hsize 10cm
\offinterlineskip\settabs\+\indent
\vtrule
\hskip 1.2cm &\vtrule 
\hskip 5.2cm &\vtrule 
\hskip 2.8cm &\vtrule 
\cr\hrule 
\+\vtrule
\hfill\,Symbol\hfill&\vtrule
\hfill Name\hfill &\vtrule
\hfill Spin-tensorial\hfill &\vtrule\cr
\vskip -0.2cm
\+\vtrule
\hfill &\vtrule
\hfill \hfill&\vtrule
\hfill type\hfill&\vtrule\cr\hrule
\+\vtrule
\hfill $\bold g$\hfill&\vtrule
\hfill Metric tensor\hfill&\vtrule
\hfill $(0,0|0,0|0,2)$\hfill&\vtrule\cr\hrule
\+\vtrule
\hfill $\bold d$\hfill&\vtrule
\hfill Skew-symmetric metric tensor\hfill&\vtrule
\hfill $(0,2|0,0|0,0)$\hfill&\vtrule\cr\hrule
\+\vtrule
\hfill$\bold H$\hfill&\vtrule
\hfill Chirality operator\hfill&\vtrule
\hfill $(1,1|0,0|0,0)$\hfill&\vtrule\cr\hrule
\+\vtrule
\hfill$\bold D$\hfill&\vtrule
\hfill Dirac form\hfill&\vtrule
\hfill $(0,1|0,1|0,0)$\hfill&\vtrule\cr\hrule
\+\vtrule
\hfill$\boldsymbol\gamma$\hfill&\vtrule
\hfill Dirac $\gamma$-field\hfill&\vtrule
\hfill $(1,1|0,0|1,0)$\hfill&\vtrule\cr\hrule
}\quad
\mytag{1.2}
$$
Note that the metric tensor $\bold g$ inherited from $M$ is interpreted
as a spin-tensorial field of $DM$ in the table \mythetag{1.2}.
\head
2. The electro-weak bundles $\uU\!M$ and $S\uuU\!M$. 
\endhead
     The first bundle $\uU\!M$ of the Standard Model is a
one-dimensional complex bundle over the space-time manifold $M$. It is
equipped with a Hermitian scalar product
$$
\hskip -2em
\uD(\bold X,\bold Y)=\uD_{11}\,\overline{X^1}\,Y^{\kern 0.2pt
\lower 1.2pt\hbox{$\ssize 1$}}=\overline{\uD(\bold Y,\bold X)}.
\mytag{2.1}
$$
(compare \mythetag{2.1} with the formula \thetag{5.18} in \mycite{2}).
Here $X^1$ and $Y^1$ are the components of two vector fields of the
bundle $\uU\!M$ in some frame $(U,\bolduPsi_1)$. Any frame
of $\uU\!M$ is a smooth section of this bundle nonzero at
each point of some open domain $U\subset M$.\par
     Tensor fields for $\uU\!M$ are defined in a way similar to those for
the Dirac bundle. In addition to $\uU\!M$, we construct its conjugate and
Hermitian conjugate bundles:
$$
\xalignat 3
&\hskip -2em
\uU^*\!M,
&&\uU^{\sssize\dagger}\!M,
&&\uU^{*\sssize\dagger}\!M=\uU^{{\sssize\dagger}*}\!M.
\mytag{2.2}
\endxalignat
$$
Then, using \mythetag{2.2}, we define the following tensor products:
$$
\gather
\hskip -2em
\uU^\varepsilon_\eta M=\overbrace{\uU\!M\otimes\ldots\otimes 
\uU\!M}^{\text{$\varepsilon$ times}}\otimes
\underbrace{\uU^*\!M\otimes\ldots\otimes 
\uU^*\!M}_{\text{$\eta$ times}}\,,
\mytag{2.3}\\
\hskip -2em
\bar{\uU}^{\kern 0.2pt\lower 2.1pt\hbox{$\ssize\sigma$}}_\zeta M
=\overbrace{\uU^{{\sssize\dagger}*}\!M\otimes
\ldots\otimes \uU^{{\sssize\dagger}*}\!M}^{\text{$\sigma$ times}}
\otimes\underbrace{\uU^{\sssize\dagger}\!M)\otimes\ldots\otimes 
\uU^{\sssize\dagger}\!M}_{\text{$\zeta$ times}}.
\mytag{2.4}
\endgather
$$     
The tensor product of the above two tensor bundles \mythetag{2.3} and
\mythetag{2.4} is denoted as
$$
\pagebreak
\hskip -2em
\uU^\varepsilon_\eta\bar{\uU}^{\kern 0.2pt\lower 2.1pt\hbox{$\ssize
\sigma$}}_\zeta M=
\uU^\varepsilon_\eta M\otimes\bar{\uU}^{\kern 0.2pt
\lower 2.1pt\hbox{$\ssize\sigma$}}_\zeta M.
\mytag{2.5}
$$
Due to \mythetag{2.5} the tensorial type of a $\uU\!M$-tensor is given
by four numbers $(\varepsilon,\eta|\sigma,\zeta)$. In particular, the
Hermitian scalar product \mythetag{2.1} is given by a tensor field
$\bolduD$ of the type $(0,1|0,1)$. This field could be upended to those 
in the table \mythetag{1.2}.\par
\mydefinition{2.1} A frame $(U,\bolduPsi_1)$ of the electro-weak bundle
$\uU\!M$ is called an {\it orthonormal frame} with respect to $\bolduD$
if $\uD_{11}=1$ in this frame.
\enddefinition
    Assume that $(U,\bolduPsi_1)$ and $(\tilde U,\tilde{\bolduPsi}_1)$
are two orthonormal frames with respect to $\bolduD$ and assume that
their domains are overlapping: $U\cap\tilde U\neq\varnothing$. Then
$$
\hskip -2em
\tilde{\bolduPsi}_1=\goth S^1_1\,\bolduPsi_1\text{, \ where \ }
|\goth S^1_1|=1.
\mytag{2.6}
$$
The coefficient $\goth S^1_1$ in \mythetag{2.6} is interpreted as a
unitary matrix $\goth S\in\MatGrU(1)$.\par
    The electro-weak bundle $\uU\!M$ is a complex bundle over the real
manifold $M$. For this reason it is equipped with the involution of 
complex conjugation $\tau$:
$$
\hskip -2em
\CD
@>\tau>>\\
\vspace{-4ex}
\uU^\varepsilon_\eta\bar{\uU}^{\kern 0.2pt\lower 2.1pt\hbox{$\ssize
\sigma$}}_\zeta M@.
\uU^\sigma_\zeta\bar{\uU}^{\kern 0.2pt\lower 2.1pt
\hbox{$\ssize\varepsilon$}}_\eta M.\\
\vspace{-4.2ex}
@<<\tau< 
\endCD
\mytag{2.7}
$$
The same is true for the Dirac bundle $DM$ (see the formula \thetag{1.2} 
in \mycite{3}).\par
    The second electro-weak bundle $S\uuU\!M$ is a little bit more
complicated. It is a two-dimensional complex bundle over the space-time
manifold $M$. Like $\uU\!M$, this bundle is equipped with a Hermitian
scalar product
$$
\hskip -2em
\uuD(\bold X,\bold Y)=\sum^2_{i=1}\sum^2_{\bar j=1}
\uuD_{i\bar j}\,\overline{X^{\bar j}}\,Y^{\kern 0.2pt
\lower 1.2pt\hbox{$\ssize i$}}=\overline{\uuD(\bold Y,\bold X)}.
\mytag{2.8}
$$
Tensor bundles associated with $S\uuU\!M$ are introduced by the following
formulas:
$$
\align
&\hskip -2em
S\uuU^\pi_\rho M=\overbrace{S\uuU\!M\otimes\ldots\otimes 
S\uuU\!M}^{\text{$\pi$ times}}\otimes
\underbrace{S\uuU^*\!M\otimes\ldots\otimes 
S\uuU^*\!M}_{\text{$\rho$ times}}\,,
\mytag{2.9}\\
&\hskip -2em
\overline{S\uuU}^{\kern 0.2pt\lower 2.1pt\hbox{$\ssize\omega$}}_\mu M
=\overbrace{S\uuU^{{\sssize\dagger}*}\!M\otimes
\ldots\otimes S\uuU^{{\sssize\dagger}*}\!M}^{\text{$\omega$ times}}
\otimes\underbrace{S\uuU^{\sssize\dagger}\!M)\otimes\ldots\otimes 
S\uuU^{\sssize\dagger}\!M}_{\text{$\mu$ times}}.
\mytag{2.10}\\
\vspace{2ex}
&\hskip -2em
S\uuU^\pi_\rho\overline{S\uuU}^{\kern 0.2pt\lower 2.1pt
\hbox{$\ssize\omega$}}_\mu M=S\uuU^\pi_\rho M\otimes
\overline{S\uuU}^{\kern 0.2pt\lower 2.1pt\hbox{$\ssize
\omega$}}_\mu M.
\mytag{2.11}
\endalign
$$
The formulas \mythetag{2.9}, \mythetag{2.10}, and \mythetag{2.11} are
analogous to \mythetag{2.3}, \mythetag{2.4}, and \mythetag{2.5}. As
we see in \mythetag{2.11}, the tensorial type of a $S\uuU\!M$-tensor 
is determined by four numbers $(\pi,\rho|\omega,\mu)$. Like \mythetag{2.1},
the Hermitian scalar product \mythetag{2.8} is determined by a tensor-field
$\bolduuD$ of the type $(0,1|0,1)$. We call it the {\it Hermitian metric 
tensor\/} of the electro-weak bundle $S\uuU\!M$.\par
     Apart from $\bolduuD$, there is another basic tensor field associated
with the bundle $S\uuU\!M$. It is a skew-symmetric tensor of the type
$(0,2|0,0)$. We denote it by $\bolduud$ and call the {\it skew-symmetric 
metric tensor}. This tensor is assumed to be nonzero at each point $p$ of
the space-time manifold $M$. Therefore one can introduce the {\it dual
skew-symmetric metric tensor}. By tradition we denote it with the same
symbol $\bolduud$. The components $\uud^{ij}$ of the dual skew-symmetric
metric tensor form the matrix inverse to the matrix formed by the
components $\uud_{ij}$ of the initial tensor $\bolduud$:
$$
\xalignat 2
&\hskip -2em
\sum^2_{j=1}\uud^{ij}\,\uud_{jk}=\delta^i_k,
&&\sum^2_{j=1}\uud_{ij}\,\uud^{jk}=\delta^k_i.
\mytag{2.12}
\endxalignat
$$
The mutually inverse matrices $\uud_{ij}$ and $\uud^{ij}$ from
\mythetag{2.12} are used in index raising and index lowering 
procedures for $S\uuU\!M$-tensors. Let's apply them to the tensor
$\bolduuD$:
$$
\hskip -2em
\uuD^{i\bar j}=\sum^2_{p=1}\sum^2_{\bar q=1}
\uud^{\kern 0.5pt ip}\,\overline{\uud^{\kern 0.4pt\raise 1.2pt
\hbox{$\ssize\bar j\bar q$}}}\,\uuD_{p\kern 0.5pt\bar q}.
\mytag{2.13}
$$
\mydefinition{2.2} The skew-symmetric metric tensor $\bolduud$ is called
{\it concordant with the Hermitian scalar product \mythetag{2.8}} if the
matrix \mythetag{2.13} is inverse to the matrix $\uuD_{i\bar j}$ in the
sense of the following equalities:
$$
\xalignat 2
&\hskip -2em
\sum^2_{\bar a=1}\uuD^{i\bar a}\,\uuD_{j\bar a}=\delta^i_j,
&&\sum^2_{a=1}\uuD_{a\bar j}\,\uuD^{a\bar i}=\delta^{\bar i}_{\bar j}.
\mytag{2.14}
\endxalignat
$$
\enddefinition
    Note that almost the same relationships are valid for the 
spin-tensorial fields $\bold d$ and $\bold D$ associated with 
the Dirac bundle (see formulas \thetag{6.24} and \thetag{6.25}
in \mycite{2}).
\mydefinition{2.3} A frame $(U,\,\bolduuPsi_1,\,\bolduuPsi_2)$ of the
bundle $S\uuU\!M$ is called an {\it orthonormal frame} with respect 
to $\bolduuD$ if $\bolduuD$ is given by the unit matrix in this frame:
$$
\hskip -2em
\uuD_{i\bar j}=\Vmatrix 1 & 0\\0 & 1\endVmatrix.
\mytag{2.15}
$$
\enddefinition
\mydefinition{2.4} A frame $(U,\,\bolduuPsi_1,\,\bolduuPsi_2)$ of the
bundle $S\uuU\!M$ is called an {\it orthonormal frame} with respect 
to $\bolduud$ if $\bolduud$ is given by the following matrix in this
frame:
$$
\hskip -2em
\uud_{ij}=\Vmatrix 0 & 1\\-1 & 0\endVmatrix.
\mytag{2.16}
$$
\enddefinition
\mydefinition{2.5} A frame $(U,\,\bolduuPsi_1,\,\bolduuPsi_2)$ of the
electro-weak bundle $S\uuU\!M$ is called an {\it orthonormal frame} if 
it is orthonormal with respect to $\bolduuD$ and $\bolduud$ simultaneously,
i\.\,e\. if both equalities \mythetag{2.15} and \mythetag{2.16} are valid
in this frame.
\enddefinition
\mytheorem{2.1} The tensor fields $\bolduuD$ and $\bolduud$ associated 
with the electro-weak bundle $S\uuU\!M$ are concordant in the sense of the
definition~\mythedefinition{2.2} if and only if for each point $p\in M$
there is an orthonormal frame $(U,\,\bolduuPsi_1,\,\bolduuPsi_2)$ in some
neighborhood $U$ of $p$.
\endproclaim
     The proof of this theorem is left to the reader. It is rather simple.
Note that the concordance of $\bold d$ and $\bold D$ in the case of the 
Dirac bundle $DM$ follows from the definition of these spin-tensorial 
fields. In the present case it is introduced as an additional requirement
for the fields $\bolduuD$ and $\bolduud$ and for the bundle $S\uuU\!M$
itself.\par
     Assume that we have two orthonormal frames $(U,\,\bolduuPsi_1,
\,\bolduuPsi_2)$ and $(\tilde U,\,\tilde{\bolduuPsi}_1,
\,\tilde{\bolduuPsi}_2)$ of the bundle $S\uuU\!M$ with overlapping domains
$U\cap\tilde U\neq\varnothing$. Then
$$
\hskip -2em
\tilde{\bolduuPsi}_i=\sum^2_{j=1}\goth S^j_i\,\bolduuPsi_j
\mytag{2.17}
$$
at each point $p\in U\cap\tilde U$. The transition matrix $\goth S$ with
the components $\goth S^j_i$ in \mythetag{2.17} is a unitary matrix
with $\det\goth S=1$, i\.\,e\. $\goth S\in\MatGrSU(2)$.\par
    Like $\uU\!M$, the second electro-weak bundle $S\uuU\!M$ is equipped
with the involution of complex conjugation. We denote it with the same
symbol $\tau$:
$$
\hskip -2em
\CD
@>\tau>>\\
\vspace{-4ex}
S\uuU^\pi_\rho\overline{S\uuU}^{\kern 0.2pt\lower 2.1pt
\hbox{$\ssize\omega$}}_\mu M@.
S\uuU^\omega_\mu\overline{S\uuU}^{\kern 0.2pt\lower 2.1pt
\hbox{$\ssize\pi$}}_\rho M.\\
\vspace{-4.2ex}
@<<\tau< 
\endCD
\mytag{2.18}
$$
The diagram \mythetag{2.18} is analogous to \mythetag{2.7} and to the 
diagram \thetag{1.2} in \mycite{3}.
\head
3. The color bundle $S\uuuU\!M$. 
\endhead
     The color bundle $S\uuuU\!M$ is used to describe the color states
of quarks (see \mycite{5}). This is a three-dimensional complex bundle
over the real space-time manifold $M$. For this reason it is equipped
with the involution of complex conjugation $\tau$:
$$
\hskip -2em
\CD
@>\tau>>\\
\vspace{-4ex}
S\uuuU^\varkappa_\lambda\overline{S\uuuU}^{\kern 0.2pt\lower 2.1pt
\hbox{$\ssize\varsigma$}}_\chi M@.
S\uuuU^\varsigma_\chi\overline{S\uuuU}^{\kern 0.2pt\lower 2.1pt
\hbox{$\ssize\varkappa$}}_\lambda M.\\
\vspace{-4.2ex}
@<<\tau< 
\endCD
\mytag{3.1}
$$
The tensor bundles $S\uuuU^\varkappa_\lambda\overline{S
\uuuU}^{\kern 0.2pt\lower 2.1pt\hbox{$\ssize\varsigma$}}_\chi M$ 
and $S\uuuU^\varsigma_\chi\overline{S\uuuU}^{\kern 0.2pt\lower 2.1pt
\hbox{$\ssize\varkappa$}}_\lambda M$ in \mythetag{3.1} are defined 
as follows:
$$
\align
&\hskip -2em
S\uuuU^\varkappa_\lambda M=\overbrace{S\uuuU\!M\otimes\ldots\otimes 
S\uuuU\!M}^{\text{$\varkappa$ times}}\otimes
\underbrace{S\uuuU^*\!M\otimes\ldots\otimes 
S\uuuU^*\!M}_{\text{$\lambda$ times}}\,,
\mytag{3.2}\\
&\hskip -2em
\overline{S\uuuU}^{\kern 0.2pt\lower 2.1pt\hbox{$\ssize\varsigma$}}_\chi M
=\overbrace{S\uuuU^{{\sssize\dagger}*}\!M\otimes
\ldots\otimes S\uuuU^{{\sssize\dagger}*}\!M}^{\text{$\varsigma$ times}}
\otimes\underbrace{S\uuuU^{\sssize\dagger}\!M)\otimes\ldots\otimes 
S\uuuU^{\sssize\dagger}\!M}_{\text{$\chi$ times}}.
\mytag{3.3}\\
\vspace{2ex}
&\hskip -2em
S\uuuU^\varkappa_\lambda\overline{S\uuuU}^{\kern 0.2pt\lower 2.1pt
\hbox{$\ssize\varsigma$}}_\chi M=S\uuuU^\varkappa_\lambda M\otimes
\overline{S\uuuU}^{\kern 0.2pt\lower 2.1pt\hbox{$\ssize
\varsigma$}}_\chi M.
\mytag{3.4}
\endalign
$$
The formulas \mythetag{3.2}, \mythetag{3.3}, \mythetag{3.4} are 
analogous to \mythetag{2.3}, \mythetag{2.4}, \mythetag{2.5} and
\mythetag{2.9}, \mythetag{2.10}, \mythetag{2.11}. In addition to 
$\tau$, the color bundle $S\uuuU\!M$ is equipped with two tensorial 
fields $\bolduuuD$ and $\bolduuud$. The tensorial field $\bolduuuD$
is introduced through the Hermitian scalar product
$$
\hskip -2em
\uuuD(\bold X,\bold Y)=\sum^2_{i=1}\sum^2_{\bar j=1}
\uuuD_{i\bar j}\,\overline{X^{\bar j}}\,Y^{\kern 0.2pt
\lower 1.2pt\hbox{$\ssize i$}}=\overline{\uuuD(\bold Y,\bold X)}.
\mytag{3.5}
$$
It is a field of the type $(0,1|0,1)$. The second tensorial field 
$\bolduuud$ is a nonzero completely skew-symmetric field of the type
$(0,3|0,0)$. For its components we have:
$$
\xalignat 3
&\uuud_{ijk}=-\uuud_{ji\kern 0.4pt k},
&&\uuud_{ijk}=-\uuud_{i\kern 0.4pt kj},
&&\uuud_{ijk}=-\uuud_{kj\kern 0.4pt i}.
\qquad
\mytag{3.6}
\endxalignat
$$
The formulas \mythetag{3.6} means that all of the components of the tensor
$\bolduuud$ can be expressed through $\uuud_{123}$. The following
definition is based on this feature.
\mydefinition{3.1} A frame $(U,\,\bolduuuPsi_1,\,\bolduuuPsi_2,
\,\bolduuuPsi_3)$ of the bundle $S\uuuU\!M$ is called an {\it orthonormal 
frame} with respect to $\bolduuud$ if $\uuud_{123}=1$ in this frame.
\enddefinition
\mydefinition{3.2} A frame $(U,\,\bolduuuPsi_1,\,\bolduuuPsi_2,
\,\bolduuuPsi_3)$ of the bundle $S\uuuU\!M$ is called an {\it orthonormal
frame} with respect to $\bolduuuD$ if $\bolduuuD$ is given by the unit 
matrix in this frame:
$$
\uuuD_{i\bar j}=\Vmatrix 1 & 0 & 0\\
\vspace{0.5ex}
0 & 1 & 0\\
\vspace{0.5ex}
0 & 0 & 1\endVmatrix.
$$
\enddefinition
Let $\bold b$ be a completely skew-symmetric tensor field of the type
$(3,0|0,0)$ associated with the color bundle $S\uuuU\!M$. Then we have
the identities similar to \mythetag{3.6}:
$$
\xalignat 3
&b^{ijk}=-b^{ji\kern 0.4pt k},
&&b^{ijk}=-b^{i\kern 0.4pt kj},
&&b^{ijk}=-b^{kj\kern 0.4pt i}.
\qquad
\mytag{3.7}
\endxalignat
$$
The identities \mythetag{3.7} mean that the components of the tensor
$\bold b$ are completely determined if the component $b^{123}$ is fixed.
\mydefinition{3.3} A completely skew-symmetric tensor $\bold b$ of the
type $(3,0|0,0)$ is called {\it inverse to the tensor $\bolduuud$} if 
\ $\uuud_{123}\cdot b^{123}=1$ \ for any frame $(U,\,\bolduuuPsi_1,
\,\bolduuuPsi_2,\,\bolduuuPsi_3)$.
\enddefinition
For the sake of economy the tensor $\bold b$ inverse to the tensor 
$\bolduuud$ is denoted by the same symbol, i\.\,e\. we write
$\bold b=\bolduuud$. This usage of symbols makes no confusion since
the initial and the inverse tensors are of different types $(0,3|0,0)$ 
and $(3,0|0,0)$ respectively.
\mydefinition{3.4} The completely skew-symmetric tensor field $\bolduuud$ 
is called {\it concordant with the Hermitian scalar product \mythetag{3.5}}
if the following identity is fulfilled:
$$
\hskip -2em
\sum^3_{i=1}\sum^3_{j=1}\sum^3_{k=1}\uuud^{\kern 0.4pt ijk}
\,\uuuD_{i\bar i}\,\uuuD_{j\bar j}\,\uuuD_{k\bar k}
=\overline{\uuud_{\kern 0.4pt\bar i\kern 0.2pt\bar j\bar k}}.
\mytag{3.8}
$$
\enddefinition
Note that the equalities \mythetag{2.13} and \mythetag{2.14} can be
transformed to the following equality for the components of the
electro-weak tensors $\bolduuD$ and $\bolduud$:
$$
\hskip -2em
\sum^2_{i=1}\sum^2_{j=1}\uud^{\kern 0.4pt ij}
\,\uuD_{i\bar i}\,\uuD_{j\bar j}
=\overline{\uud_{\kern 0.4pt\bar i\kern 0.2pt\bar j}}.
\mytag{3.9}
$$
The equality \mythetag{3.9} is a two-dimensional analog of the 
equality \mythetag{3.8}. Therefore, the definition~\mythedefinition{2.2}
also is a two-dimensional analog of the definition~\mythedefinition{3.4}.
\mydefinition{3.5} A frame $(U,\,\bolduuuPsi_1,\,\bolduuuPsi_2,
\,\bolduuuPsi_3)$ of the color bundle $S\uuuU\!M$ is called an 
{\it orthonormal frame} if it is orthonormal with respect to 
$\bolduuuD$ and $\bolduuud$ simultaneously.
\enddefinition
\mytheorem{3.1} The tensor fields $\bolduuuD$ and $\bolduuud$ associated 
with the color bundle $S\uuuU\!M$ are concordant in the sense of the
definition~\mythedefinition{3.4} if and only if for each point $p\in M$
there is an orthonormal frame $(U,\,\bolduuuPsi_1,\,\bolduuuPsi_2,
\,\bolduuuPsi_3)$ in some neighborhood $U$ of $p$.
\endproclaim
\noindent This theorem~\mythetheorem{3.1} is similar to the
theorem~\mythetheorem{2.1}. Its proof is also left to the reader.
\par
    Assume that we have two orthonormal frames of the color bundle 
$S\uuuU\!M$. Let's denote them $(U,\,\bolduuuPsi_1,\,\bolduuuPsi_2,
\,\bolduuuPsi_3)$ and $(\tilde U,\,\tilde {\bolduuuPsi}_1,\,
\tilde{\bolduuuPsi}_2,\,\tilde{\bolduuuPsi}_3)$. Assume that their
domains $U$ and $\tilde U$ are overlapping, i\.\,e\. $U\cap\tilde U
\neq\varnothing$. Then at each point $p\in U\cap\tilde U$ we have
the transition formula relating two frames $(U,\,\bolduuuPsi_1,\,
\bolduuuPsi_2,\,\bolduuuPsi_3)$ and $(\tilde U,\,\tilde{\bolduuuPsi}_1,
\,\tilde{\bolduuuPsi}_2,\,\tilde{\bolduuuPsi}_3)$:
$$
\hskip -2em
\tilde{\bolduuuPsi}_i=\sum^3_{j=1}\goth S^j_i\,\bolduuuPsi_j.
\mytag{3.10}
$$
The transition matrix $\goth S$ with the components $\goth S^j_i$ in
\mythetag{3.10} is a unitary matrix with $\det\goth S=1$, i\.\,e\. it
is an element of the special unitary group: $\goth S\in\MatGrSU(3)$.
\par
{\bf Conclusion}. As we see now, the structural groups for all of the 
above bundles are revealed through the transition matrices relating 
orthonormal frames, while orthonormal frames themselves are defined with
the use of some special tensorial fields associated with these bundles.
\par
\head
4. Leptons and quarks.
\endhead
    Basic elementary particles in the Standard Model are subdivided
into two groups, the first group is formed by leptons and the second
group is formed by quarks. Other particles arise as quanta of gauge 
fields. Leptons are described by the $\MatGrU(1)\times
\MatGrSU(2)$ symmetry. Their wave functions are associated with the 
following tensor bundles:
$$
\hskip -2em
S\uuU^\pi_\rho\overline{S\uuU}^{\kern 0.2pt\lower 2.1pt
\hbox{$\ssize\omega$}}_\mu M\otimes
\uU^\varepsilon_\eta\bar{\uU}^{\kern 0.2pt\lower 
2.1pt\hbox{$\ssize\sigma$}}_\zeta M\otimes
D^\alpha_\beta\bar D^\nu_\gamma T^m_n M.
\mytag{4.1}
$$
Smooth sections of the bundle \mythetag{4.1} are called smooth
spin-tensorial fields. Apart from the spin-tensorial type 
$(\alpha,\beta|\nu,\gamma|m,n)$, they are characterized by the
electro-weak type $(\pi,\rho|\omega\mu|\varepsilon,\eta|\sigma\zeta)$.
The following table is an electro-weak addendum to the table
\mythetag{1.2} containing the basic electro-weak fields:
$$
\vcenter{\hsize 10cm
\offinterlineskip\settabs\+\indent
\vtrule
\hskip 1.2cm &\vtrule 
\hskip 5.2cm &\vtrule 
\hskip 3.4cm &\vtrule 
\cr\hrule 
\+\vtrule
\hfill\,Symbol\hfill&\vtrule
\hfill Name\hfill &\vtrule
\hfill Spin-tensorial and\hfill &\vtrule\cr
\vskip -0.2cm
\+\vtrule
\hfill &\vtrule
\hfill \hfill&\vtrule
\hfill electro-weak types\hfill&\vtrule\cr\hrule
\+\vtrule
\hfill $\bolduD$\hfill&\vtrule
\hfill Hermitian metric tensor\hfill&\vtrule
\hfill $(0,0|0,0|0,0)$\hfill&\vtrule\cr
\+\vtrule
\hfill &\vtrule\hfill&\vtrule
\hfill $(0,0|0,0|0,1|0,1)$\hfill&\vtrule\cr\hrule
\+\vtrule
\hfill$\bolduuD$\hfill&\vtrule
\hfill Hermitian metric tensor\hfill&\vtrule
\hfill $(0,0|0,0|0,0)$\hfill&\vtrule\cr
\+\vtrule
\hfill &\vtrule\hfill&\vtrule
\hfill $(0,1|0,1|0,0|0,0)$\hfill&\vtrule\cr\hrule
\+\vtrule
\hfill$\bolduud$\hfill&\vtrule
\hfill Skew-symmetric metric tensor\hfill&\vtrule
\hfill $(0,0|0,0|0,0)$\hfill&\vtrule\cr
\+\vtrule
\hfill &\vtrule\hfill&\vtrule
\hfill $(0,2|0,0|0,0|0,0)$\hfill&\vtrule\cr\hrule
}\quad
\mytag{4.2}
$$\par
     Quarks are described by the $\MatGrU(1)\times\MatGrSU(2)\times 
\MatGrSU(3)$ symmetry. Their wave functions are associated with the 
following tensor bundles:
$$
\hskip -2em
S\uuuU^\varkappa_\lambda\overline{S\uuuU}^{\kern 0.2pt\lower 2.1pt
\hbox{$\ssize\varsigma$}}_\chi M\otimes
S\uuU^\pi_\rho\overline{S\uuU}^{\kern 0.2pt\lower 2.1pt
\hbox{$\ssize\omega$}}_\mu M\otimes
\uU^\varepsilon_\eta\bar{\uU}^{\kern 0.2pt\lower 2.1pt
\hbox{$\ssize\sigma$}}_\zeta M\otimes
D^\alpha_\beta\bar D^\nu_\gamma T^m_n M.
\mytag{4.3}
$$
Like in the case of the bundle \mythetag{4.1}, smooth sections of the
bundle \mythetag{4.3} are called smooth spin-tensorial fields. Apart 
from the spin-tensorial type $(\alpha,\beta|\nu,\gamma|m,n)$, they are
characterized by the color type $(\varkappa,\lambda|\varsigma,\chi|
\pi,\rho|\omega,\mu|\varepsilon,\eta|\sigma,\zeta)$. The following table 
is a color addendum to the tables \mythetag{1.2} and
\mythetag{4.2}:
$$
\vcenter{\hsize 10cm
\offinterlineskip\settabs\+\indent
\vtrule
\hskip 1.2cm &\vtrule 
\hskip 4.2cm &\vtrule 
\hskip 4.4cm &\vtrule 
\cr\hrule 
\+\vtrule
\hfill\,Symbol\hfill&\vtrule
\hfill Name\hfill &\vtrule
\hfill Spin-tensorial and\hfill &\vtrule\cr
\vskip -0.2cm
\+\vtrule
\hfill &\vtrule
\hfill \hfill&\vtrule
\hfill color types\hfill&\vtrule\cr\hrule
\+\vtrule
\hfill $\bolduD$\hfill&\vtrule
\hfill Hermitian metric tensor\hfill&\vtrule
\hfill $(0,0|0,0|0,0)$\hfill&\vtrule\cr
\+\vtrule
\hfill &\vtrule\hfill&\vtrule
\hfill $(0,0|0,0|0,0|0,0|0,1|0,1)$\hfill&\vtrule\cr\hrule
\+\vtrule
\hfill$\bolduuD$\hfill&\vtrule
\hfill Hermitian metric tensor\hfill&\vtrule
\hfill $(0,0|0,0|0,0)$\hfill&\vtrule\cr
\+\vtrule
\hfill &\vtrule\hfill&\vtrule
\hfill $(0,0|0,0|0,1|0,1|0,0|0,0)$\hfill&\vtrule\cr\hrule
\+\vtrule
\hfill$\bolduud$\hfill&\vtrule
\hfill Skew-symmetric\hfill&\vtrule
\hfill $(0,0|0,0|0,0)$\hfill&\vtrule\cr
\+\vtrule
\hfill &\vtrule\hfill metric tensor\hfill&\vtrule
\hfill $(0,0|0,0|0,2|0,0|0,0|0,0)$\hfill&\vtrule\cr\hrule
\+\vtrule
\hfill$\bolduuuD$\hfill&\vtrule
\hfill Hermitian metric tensor\hfill&\vtrule
\hfill $(0,0|0,0|0,0)$\hfill&\vtrule\cr
\+\vtrule
\hfill &\vtrule\hfill&\vtrule
\hfill $(0,1|0,1|0,0|0,0|0,0|0,0)$\hfill&\vtrule\cr\hrule
\+\vtrule
\hfill$\bolduuud$\hfill&\vtrule
\hfill Completely\hfill&\vtrule
\hfill $(0,0|0,0|0,0)$\hfill&\vtrule\cr
\+\vtrule
\hfill &\vtrule\hfill skew-symmetric tensor\hfill&\vtrule
\hfill $(0,3|0,0|0,0|0,0|0,0|0,0)$\hfill&\vtrule\cr\hrule
}\quad
\mytag{4.4}
$$\par
    Note that $\MatGrU(1)$ and $\MatGrSU(2)$-bundles in the case of
leptons and in the case of quarks could be different (non-isomorphic)
bundles. In this case the field $\bolduD$, $\bolduuD$, $\bolduud$ in
the table \mythetag{4.4} are different from $\bolduD$, $\bolduuD$,
$\bolduud$ in the previous table \mythetag{4.2}. But even if it is so, 
it is convenient to denote these different fields by the same symbols.

\enddocument
\end